\documentclass[12pt]{article}
\usepackage{amssymb,amsmath}
\usepackage{cases}
\usepackage{amsfonts}
\usepackage{color,xcolor}
\usepackage[left=2.0cm,right=2.0cm,top=2.0cm,bottom=2.0cm]{geometry}
\usepackage[colorlinks,citecolor=blue,urlcolor=blue]{hyperref}

\newtheorem{theorem}{Theorem}[section]

\newtheorem{lemma}{Lemma}[section]

\newtheorem{remark}{Remark}[section]

\newcommand{\bal}{\begin{align}}
\newcommand{\bbal}{\begin{align*}}
\newcommand{\beq}{\begin{equation}}
\newcommand{\eeq}{\end{equation}}
\newcommand{\bca}{\begin{cases}}
\newcommand{\eca}{\end{cases}}
\def\div{\mathord{{\rm div}}}
\newcommand{\pa}{\partial}
\newcommand{\fr}{\frac}
\newcommand{\na}{\nabla}

\newcommand{\cd}{\cdot}
\newcommand{\ep}{\varepsilon}
\newcommand{\dd}{\mathrm{d}}

\newcommand{\R}{\mathbb{R}}

\newcommand{\les}{\lesssim}

\newcommand{\bi}{\Big}

\begin{document}

\title{Global smooth solutions of the generalized MHD equations with large data}

\author{Jinlu Li$^{1}$\footnote{E-mail: lijinlu@gnnu.cn}\quad  and Yanghai Yu$^{2}$\footnote{E-mail: yuyanghai214@sina.com( Corresponding author)}\\
\small $^1$\it School of Mathematics and Computer Sciences, Gannan Normal University, Ganzhou 341000, China\\
\small $^2$\it School of Mathematics and Statistics, Anhui Normal University, Wuhu, Anhui, 241002, China}

\date{}

\maketitle\noindent{\hrulefill}

{\bf Abstract:} In this paper, we consider the Cauchy problem of the multi-dimensional generalized MHD system in the whole space and construct global smooth solutions with a class of large initial data by exploring the structure of the nonlinear term. Precisely speaking, our choice of special initial data whose $L^{\infty}$ norm can be arbitrarily large allows to generate global-in-time solutions to the generalized MHD system.


\vskip0mm\noindent{\hrulefill}

\section{Introduction}\label{sec1}
This paper focuses on the following generalized incompressible magnetohydrodynamics (MHD) equations
\begin{eqnarray}\label{g-mhd}
        \left\{\begin{array}{ll}
          \partial_tu+u\cd\na u+\mu\Lambda^{\alpha} u+\na p=b\cd\na b,& x\in \R^d,t>0,\\
          \partial_tb+u\cd\na b+\nu\Lambda^{\beta} b=b\cd\na u,& x\in \R^d,t>0,\\
         \div u=\div b=0,& x\in \R^d,t\geq0,\\
          (u,b)|_{t=0}=(u_0,b_0),& x\in \R^d,\end{array}\right.
        \end{eqnarray}
where $u$ and $b$ denote the divergence free velocity field and magnetic field, respectively, $p\in \R$ is the scalar pressure. The parameter $\mu$ denotes the kinematic viscosity coefficient of the fluid and $\nu$ is the magnetic diffusivity coefficient. The fractional power operator $\Lambda^{\gamma}$ with $0<\gamma<2$ is defined by Fourier multiplier with symbol $|\xi|^{\gamma}$ (see e.g. \cite{Jacob 2005,Wu 2017})
\begin{eqnarray*}
  \Lambda^{\gamma}u(x)=\mathcal{F}^{-1}|\xi|^{\gamma}\mathcal{F}u(\xi).
\end{eqnarray*}
Throughout this paper, to simplify the presntation, we make the convention that by $\gamma= 0$ and $\gamma=2$ we mean that $\Lambda^{\gamma}u$ are damping term $u$ and Laplacian term $-\Delta u$, respectively. Roughly speaking, the MHD equations include the Navier-Stokes (NS) (or Euler when $\nu=0$) system as a special case, which govern the motion of electrically conducting fluids such as plasmas, liquid metals and electrolytes, and play a fundamental role in geophysics, astrophysics, cosmology and engineering (see e.g.\cite{Priest 2000,Davidson 2001,Li 2017}). From a mathematical view, the global regularity or finite time singularity for strong solutions of the 3D MHD system with large initial data is still a challenging open problem just like the 3D NS equations. Due to the profound physical background and important mathematical significance, the MHD equations have attracted considerable attention recently from the community of mathematical fluids. Let us review some important works on the MHD equations \eqref{g-mhd} which are more closely to our problem. Sermange and Temam \cite{Sermange} established the global smooth solutions to the 2D MHD equations \eqref{g-mhd} with $\alpha=\beta=2$. In the completely inviscid case ($\mu=\nu=0$), the question of whether smooth solution of the MHD equations \eqref{g-mhd} with large initial data even in $\R^2$ develops singularity in finite time remains completely open. Besides these the two extreme cases, many intermediate cases, for example, the 2D MHD equations with partial dissipation, has been studied by various authors. Fan et al.\cite{Fan 2014} solved the issue of the global regularity for the MHD equations \eqref{g-mhd} with $\mu>0,\nu>0,\alpha>0,\beta=1$ . Recently, Yuan--Zhao \cite{Yuan 2018} considered the MHD equations \eqref{g-mhd} with the dissipative operators weaker than any power of the fractional Laplacian and obtained the global regularity of the corresponding system. On the other hand, Cao et al.\cite{Cao 2014}, Jiu--Zhao \cite{Jiu 2015} established the global regularity of smooth solutions to the MHD equations \eqref{g-mhd} with $\mu=0,\nu>0,\beta>1$ by different methods. Subsequently, Agelas \cite{Agelas 2016} improved this work with the diffusion $(-\Delta)^{\beta} b (\beta > 1)$ replaced by $(-\Delta)\log^{\kappa}(e-\Delta) b (\kappa > 1)$.

Since there is no global well-posedness theory for general initial data, many literatures have been devoted to the study of global existence of smooth solutions to \eqref{g-mhd} under some smallness condition. Lin--Zhang \cite{Linf 2014} established the
global well-posedness of 3D incompressible MHD type system with initial data close to some non-trivial steady state. Similar results were also obtained in two space dimensions by Lin--Xu--Zhang \cite{Linf 2015}. Zhang \cite{Zhangt} and Ren et al.\cite{Ren 2014} provided two simplified proofs of the main result in \cite{Linf 2015}. When $\alpha=\beta=0$, Wu et al. \cite{Wu 2015} obtained that the d-dimensional MHD equations \eqref{g-mhd} always possesses a unique global solution provided that the initial datum is sufficiently small in the nonhomogeneous functional setting $H^s$ with $s>1+\fr d2$. It is also worth to mention that when $b=0$, the system \eqref{g-mhd} is reduced to the NS equations. Lei--Lin--Zhou \cite{Lei 2015} constructed a family of finite energy smooth large solutions to the NS equations with the initial data close to a Beltrami flow. Subsequently, Lin--Zhang--Zhou \cite{Liny} generalized the result in \cite{Lei 2015} for the 3-D incompressible NS equations to the case of MHD with large velocity fields and large magnetic fields. Li--Yang--Yu \cite{Lj 2019} established a class global large solution to the 2D MHD equations with damp terms whose initial energy can be arbitrarily large. However, there are few results of global well-possedness for MHD system \eqref{g-mhd} with any $\alpha>0$ and $\beta>0$. Inspired by the ideas that used in \cite{Lei 2015,Lj 2019}, we aim to prove that the system \eqref{g-mhd} with $0\leq \alpha,\beta\leq 2$ can generate unique global solutions for some class of large initial data. Our main result is stated as follows.
\begin{theorem}\label{the1.1} Let $d=2,3$ and $\ep\ll 1$. Assume that the initial data fulfills ${\rm{div}}v_0={\rm{div}}c_0={\rm{div}}U_0=0$ and
$$u_0=U_0+v_0\quad \mbox{and}\quad b_0=U_0+c_0$$
where
\begin{eqnarray}\label{Equ1.2}
\mathrm{supp} \ \hat{U}_0(\xi)\subset\mathcal{C}\triangleq\Big\{\xi \big| \ 1-\ep\leq  |\xi|\leq 1+\ep\Big\}, \qquad d=2,
\end{eqnarray}
or
\begin{eqnarray}\label{Equ1.22}
\mathrm{supp} \ \hat{U}_0(\xi)\subset\mathcal{C}\triangleq\Big\{\xi \big| \ 1-\ep\leq  |\xi|\leq 1+\ep\Big\}, \quad \Lambda U_0=\na\times U_0, \qquad d=3.
\end{eqnarray}
There exists a sufficiently small positive constant $\delta$, and a universal constant $C$ such that if
\begin{align}\label{condition}
\Big(||v_0||^2_{H^3}+||c_0||^2_{H^3}+\ep||U_0||_{L^2}(1+||\hat{U}_0||_{L^1})\Big) \exp\Big( C(||\hat{U}_0||_{L^1}+\ep||U_0||_{L^2}(1+||\hat{U}_0||_{L^1}))\Big)\leq \delta,
\end{align}
then the system \eqref{g-mhd} has a unique global solution.
\end{theorem}
\begin{remark}\label{rem1.1}
For $d=2$, we set $v_0=c_0=0$ and $U_0=(\pa_2a_0,-\pa_1a_0)^{\mathsf{T}}$ with $a_0=\ep^{-1}(\log\log\frac1\ep)^{\frac12} \chi$, where the smooth function $\hat{\chi}\in[0,1]$ satisfying
\begin{align*}
\mathrm{supp} \hat{\chi}\in \mathcal{{C}}\quad\mbox{and} \quad \hat{\chi}(\xi)=1 \quad\mbox{for} \quad \xi\in\bi[1-\frac\ep2,1+\frac\ep2\bi].
\end{align*}
Then, direct calculations show that the left side of \eqref{condition} becomes
\begin{align*}
C\ep^{\frac12}\Big(\log\log \frac1\ep\Big)^2\exp\Big(C\log\log \frac1\ep\Big).
\end{align*}
Therefore, choosing $\ep$ small enough, we deduce that the system \eqref{g-mhd} has a global solution.

Moreover, we also have
\begin{align*}
||U_0||_{L^\infty}\gtrsim ||\pa_1U_{0,2}-\pa_2U_{0,1}||_{L^\infty}=||\Delta a_0||_{L^\infty}\approx |||\xi|^2\hat{a}_0||_{L^1}\gtrsim \Big(\log\log \frac1\ep\Big)^\frac12.
\end{align*}
\end{remark}
\begin{remark}\label{rem1.2}
For $d=3$, we set $v_0=c_0=0$ and $U_0=V_0+\Lambda^{-1} \na\times V_0$ with
\begin{eqnarray*}
&V_0=\ep^{-1}\Big(\log\log\frac1\ep\Big)^\frac12\na\times
\begin{pmatrix}
a_0 \\ 0 \\ 0
\end{pmatrix}
=\ep^{-1}\Big(\log\log\frac1\ep\Big)^\frac12
\begin{pmatrix}
0 \\ \pa_3a_0 \\ -\pa_2a_0
\end{pmatrix},
\end{eqnarray*}
where the smooth function $\hat{a}_0\in[0,1]$ satisfying
\begin{align*}
\mathrm{supp} \hat{a}_0\in \mathcal{{C}}\quad\mbox{and} \quad \hat{a}_0(\xi)=1 \quad\mbox{for} \quad \xi\in\bi[1-\frac\ep2,1+\frac\ep2\bi].
\end{align*}
Here, we can show that $\mathrm{div} U_0=0$ and $\na \times U_0=\Lambda U_0$.

Moreover, we also have
\begin{eqnarray*}
\hat{U}_0=\ep^{-1}\Big(\log\log\frac1\ep\Big)^\frac12
\begin{pmatrix}
\xi^2_2+\xi^2_3 \\ -\xi_1\xi_2+i\xi_3|\xi| \\ -\xi_1\xi_3-i\xi_2|\xi|
\end{pmatrix}\frac{\hat{a}_0(\xi)}{|\xi|}.
\end{eqnarray*}
Then, direct calculations show that
\begin{align*}
||\hat{U}_0||_{L^1}\approx \Big(\log\log\frac1\ep\Big)^\frac12\quad\mbox{and}\quad||{U}_0||_{L^2}\approx \ep^{-\fr12}\Big(\log\log\frac1\ep\Big)^\frac12.
\end{align*}
Thus, the left side of \eqref{condition} becomes
\begin{align*}
C\ep^{\frac12}\Big(\log\log\frac1\ep\Big)^2\exp\Big(C\log\log \frac1\ep\Big).
\end{align*}
Therefore, choosing $\ep$ small enough, we deduce that the system \eqref{g-mhd} has a global solution.

Notice that $U_{0,1}=-\ep^{-1}\Big(\log\log\frac1\ep\Big)^\frac12(\pa^2_{2}+\pa^2_{3})\Lambda^{-1} a_0$ and $\hat{U}_{0,1}\geq 0$, we can deduce that
\bbal
||U_0||_{L^\infty}\gtrsim||U_{0,1}||_{L^\infty}\approx ||\hat{U}_{0,1}||_{L^1}\gtrsim \Big(\log\log\frac1\ep\Big)^\frac12.
\end{align*}
\end{remark}

\section{Reformulation of the System}\label{sec2}
\setcounter{equation}{0}
Let $(U,B)=(e^{-\mu t}U_0,e^{-\nu t}U_0)$ be the solutions of the following system
\begin{eqnarray}\label{app-mhd}
        \left\{\begin{array}{ll}
          \pa_tU+\mu U=0,\\
          \pa_tB+\nu B=0,\\
          (U,B)|_{t=0}=(U_0,U_{0}).\end{array}\right.
        \end{eqnarray}
Hence, \eqref{app-mhd} is equivalent to the following new system
\begin{eqnarray}\label{app1-mhd}
        \left\{\begin{array}{ll}
          \pa_t U+\mu\Lambda^\alpha U=\mu(\Lambda^\alpha-1)U\triangleq f,\\
          \pa_t B+\nu\Lambda^\beta B=\nu(\Lambda^\beta-1)B\triangleq h,\\
          \div U=\div B=0,\\
          (U,B)|_{t=0}=(U_0,U_{0}).\end{array}\right.
        \end{eqnarray}
Note that $U=e^{-\mu t}U_0$ and $B=e^{-\nu t}U_0$, we can verify easily that
\bbal
-U\cd\na B+B\cd\na U=0
\end{align*}
Introducing the new quantities $$v=u-U\quad\mbox{and}\quad c=b-B,$$ from \eqref{g-mhd}, the system \eqref{app1-mhd} can be written as follows
\begin{eqnarray}\label{app2-mhd}
        \left\{\begin{array}{ll}
\partial_tv+v\cd\na v-c\cd\na c+\mu\Lambda^\alpha v+\na p=B\cd\na c+c\cd\na B-v\cd\na U-U\cd\na v+g+f,\\
\partial_tc+v\cd\na c-c\cd\na v+\nu\Lambda^\beta c=B\cd\na v+c\cd\na U-v\cd\na B-U\cd\na c+h,\\
\div v=\div c=0,\\
(v,c)|_{t=0}=(v_0,c_0).\end{array}\right.
\end{eqnarray}
where we denote
\bbal
g=-U\cd\na U+B\cd\na B.
\end{align*}
For the case $d=2$. Notice that $\na\times U=\pa_1U_2-\pa_2U_1$, we have
\bbal
U\cd\na U_1&=U_1\pa_1U_1+U_2\pa_2U_1\\
&=\pa_1\bi(\frac{|U|^2}{2}\bi)+U_2(\pa_2U_1-\pa_1U_2)\\
&=\pa_1\bi(\frac{|U|^2+(\Delta^{-1}\na\times U)^2}{2}\bi)-U_2(1+\Delta^{-1})(\na\times U)
\end{align*}
and
\bbal
U\cd\na U_2&=\pa_2\bi(\frac{|U|^2+(\Delta^{-1}\na\times U)^2}{2}\bi)+U_1(1+\Delta^{-1})(\na\times U),
\end{align*}
which implies
\bbal
U\cd\na U=\na \bi(\frac{|U|^2+(\Delta^{-1}\na\times U)^2}{2}\bi)+U^{\perp}(1+\Delta^{-1})(\na\times U).
\end{align*}
Similarly, we also have
\bbal
B\cd\na B=\na \bi(\frac{|B|^2+(\Delta^{-1}\na\times B)^2}{2}\bi)+B^{\perp}(1+\Delta^{-1})(\na\times B).
\end{align*}
Then, we can rewrite the term $g$ as follows
\bbal
g=\na \tilde{p}+G,
\end{align*}
where
\bbal
\tilde{p}=\frac{1}{2}\bi(-|U|^2-(\Delta^{-1}\na\times U)^2+|B|^2+(\Delta^{-1}\na\times B)^2\bi),
\end{align*}
and
\bbal
G=-U^{\perp}(1+\Delta^{-1})(\na\times U)+B^{\perp}(1+\Delta^{-1})(\na\times B).
\end{align*}
For the case $d=3$. Then, we deduce from the fact $\div U=\div B=0$ that
\bbal
U\cd\na U=(\na \times U)\times U+\na\Big(\frac{|U|^2}{2}\Big)=\na\Big(\frac{|U|^2}{2}\Big)+(\Lambda U-U)\times U,\\
B\cd\na B=(\na \times B)\times B+\na\Big(\frac{|B|^2}{2}\Big)=\na\Big(\frac{|B|^2}{2}\Big)+(\Lambda B-B)\times B,
\end{align*}
which means that
\bbal
g=\na\Big(\frac{|B|^2-|U|^2}{2}\Big)+(\Lambda B-B)\times B-(\Lambda U-U)\times U.
\end{align*}
\section{The Proof of Theorem \ref{the1.1}}\label{sec3}
\setcounter{equation}{0}
Firstly, we present some estimates which will be used frequently in the proof of Theorem \ref{the1.1}.
\begin{lemma}\label{le1} \cite{Majda 2001} (Commutator estimates)
There hold that
\bbal
&\sum_{0<|\alpha|\leq 3}||[D^{\alpha},\mathbf{g}]\mathbf{f}||_{L^2}\leq C(||\mathbf{f}||_{{H}^{2}}||\na \mathbf{g}||_{L^\infty}+||\mathbf{f}||_{L^\infty}||\mathbf{g}||_{{H}^3}),
\\&\sum_{0<|\alpha|\leq 3}||[D^{\alpha},\mathbf{g}]\mathbf{f}||_{L^2}\leq C(||\na \mathbf{g}||_{L^\infty}+||\na^3 \mathbf{g}||_{L^\infty})||\mathbf{f}||_{H^2}.
\end{align*}
\end{lemma}
\begin{lemma}\label{le2} \cite{Majda 2001} (Product estimates)
For $m\in \mathbb{Z}^+$ and $m\geq2$, we have
\bbal
&\sum_{|\alpha|\leq m}||D^{\alpha}(\mathbf{f}\mathbf{g})||_{L^2}\leq C||\mathbf{f}||_{H^m}||\mathbf{g}||_{H^m},
\\&\sum_{|\alpha|\leq m}||D^{\alpha}(\mathbf{f}\mathbf{g})||_{L^2}\leq C(||\mathbf{f}||_{L^\infty}+||\na^m\mathbf{f}||_{L^\infty})||\mathbf{g}||_{H^m}.
\end{align*}
\end{lemma}
\begin{lemma}\label{lem3.1} Under the assumptions of Theorem \ref{the1.1}, it holds that
\bal\label{estimate-l}
||f||_{H^3}+||h||_{H^3}\leq Ce^{-\min\{\mu,\nu\}t}\ep||U_0||_{L^2}
\end{align}
and
\bal\label{estimate-2}
||G||_{H^3}\leq Ce^{-2\min\{\mu,\nu\}t}\ep||U_0||_{L^2}||\hat{U}_0||_{L^1}.
\end{align}
\end{lemma}
{\bf Proof of Lemma \ref{lem3.1}}\quad We just consider the two-dimensional case, since the tri-dimensional case can be dealt with in a similar manner.

For the term $f$, due to the conditions $\mathrm{supp}\ \hat{U}_0(\xi)\subset \mathcal{C}$, we can show that
\bbal
||f||^2_{H^3}=e^{-2\mu t}\int_{\mathcal{C}}(1+|\xi|^2)^3\big||\xi|^\alpha-1\big|^2|\hat{U}_0|^2\dd \xi\leq Ce^{-2\mu t}\ep^2\alpha^2||U_0||^2_{L^2}.
\end{align*}
Similar argument as the term $h$, we also have
\bbal
||h||^2_{H^3}\leq  Ce^{-\nu t}\ep^2\beta^2||U_0||^2_{L^2}.
\end{align*}
By Leibniz's formula and H\"{o}lder's inequality, then we have
\bal\label{estimate-3}
||U^{\perp}(1+\Delta^{-1})(\na\times U)||_{H^3}\lesssim& ||U||_{L^\infty}||(1+\Delta^{-1})(\na\times U)||_{H^{3}}\nonumber\\&+||U||_{H^3}||(1+\Delta^{-1})(\na\times U)||_{L^\infty}\nonumber\\
\les& e^{-2\mu t}\ep||U_0||_{L^2}||\hat{U}_0||_{L^1}.
\end{align}
Similarly, one also has
\bal\label{estimate-4}
||B^{\perp}(1+\Delta^{-1})(\na\times B)||_{H^3}&\les e^{-2\nu t}\ep||U_0||_{L^2}||\hat{U}_0||_{L^1}.
\end{align}
Combining \eqref{estimate-3} and \eqref{estimate-4} gives the desired result \eqref{estimate-2}.
We complete the proof of Lemma \ref{lem3.1}. $\Box$

{\bf Proof of Theorem \ref{the1.1}}\quad Applying $D^\ell$ on $\eqref{app2-mhd}_1$ and $\eqref{app2-mhd}_2$ respectively and taking the scalar product of them with $D^\ell v$ and $D^\ell c$, respectively, adding them together and then summing the resulting over $|\alpha|\leq 3$, we get
\bal\label{z0}
\fr12\frac{\dd}{\dd t}\Big(||v||^2_{H^3}+||c||^2_{H^3}\Big)+||\Lambda^\frac{\alpha}{2} v||^2_{H^3}+||\Lambda^\frac{\beta}{2} c||^2_{H^3}\triangleq\sum^{8}_{i=1}I_i,
\end{align}
where
\bbal
&I_1=-\sum_{0<|\ell|\leq 3}\int_{\R^d}[D^{\ell},v\cd] \na v\cd D^\ell v\dd x
-\sum_{0<|\ell|\leq 3}\int_{\R^d}[D^{\ell},v\cd] \na c\cd D^\ell c\dd x,
\\&I_2=\sum_{0<|\ell|\leq 3}\int_{\R^d}[D^{\ell},c\cd] \na c\cd D^\ell v\dd x
+\sum_{0<|\ell|\leq 3}\int_{\R^d}[D^{\ell},c\cd] \na v\cd D^\ell c\dd x,
\\&I_3=-\sum_{0<|\ell|\leq 3}\int_{\R^d}D^{\ell}(U\cd \na v)\cd D^\ell v\dd x-\sum_{0<|\ell|\leq 3}\int_{\R^d}D^{\ell}(U\cd \na c)\cd D^\ell c\dd x,
\\&I_4=\sum_{0<|\ell|\leq 3}\int_{\R^d}D^{\ell}(B\cd \na c)\cd D^{\ell}v\dd x+\sum_{0<|\ell|\leq 3}\int_{\R^d}D^{\ell}(B\cd \na v)\cd D^{\ell}c\dd x,
\\&I_5=\sum_{0\leq|\ell|\leq 3}\int_{\R^d}D^{\ell}(c\cd \na B)\cd D^{\ell}v\dd x-\sum_{0\leq|\ell|\leq 3}\int_{\R^d}D^{\ell}(v\cd \na B)\cd D^{\ell}c\dd x,
\\&I_6=\sum_{0\leq|\ell|\leq 3}\int_{\R^d}D^{\ell}(c\cd \na U)\cd D^{\ell}c\dd x-\sum_{0\leq|\ell|\leq 3}\int_{\R^d}D^{\ell}(v\cd \na U)\cd D^{\ell}v\dd x,
\\&I_7=\sum_{0\leq |\ell|\leq 3}\int_{\R^d}D^{\ell}(f+G)\cd D^{\ell} v\dd x,
\quad I_8=\sum_{0\leq |\ell|\leq 3}\int_{\R^d}D^{\ell}h\cd D^{\ell} c\dd x.
\end{align*}
Next, we need to estimate the above terms one by one.

According to the commutate estimate (See Lemma \ref{le1}), we obtain
\bal
I_1\leq&~\sum_{0<|\ell|\leq 3}||[D^{\ell},v\cd] \na v||_{L^2}||\na v||_{H^2}+\sum_{0<|\ell|\leq 3}|||[D^{\ell},v\cd]\na c||_{L^2}||\na c||_{H^2}\nonumber\\
\leq&~C||\na v||_{L^\infty}||v||_{H^3}||\na v||_{H^2}+C||v||_{H^3}||\na c||^2_{H^2}\nonumber\\
\leq&~C||v||_{H^3}\Big(||\na v||^2_{H^2}+||\na c||^2_{H^2}\Big),\label{z1}\\
I_2\leq&~\sum_{0<|\ell|\leq 3}||[D^{\ell},c\cd] \na c||_{L^2}||\na v||_{H^2}+\sum_{0<|\ell|\leq 3}|||[D^{\ell},c\cd] \na v||_{L^2}||\na c||_{H^2}\nonumber\\
\leq&~C||\na c||_{H^2}||\na v||_{H^2}||c||_{H^3}\nonumber\\
\leq&~C||c||_{H^3}\Big(||\na v||^2_{H^2}+||\na c||^2_{H^2}\Big).\label{z2}
\end{align}
Invoking the calculus inequality (See Lemma \ref{le2}), we obtain
\bal
I_3\leq&\sum_{0<|\ell|\leq 3}||[D^{\ell},U\cd] \na v||_{L^2}||\na v||_{H^2}+\sum_{0<|\ell|\leq 3}||[D^{\ell},U\cd] \na c||_{L^2}||\na c||_{H^2}\nonumber\\
\leq&C\Big(||\na U||_{L^\infty}+||\na^3 U||_{L^\infty}\Big)\Big(||v||^2_{H^3}+||c||^2_{H^3}\Big),\label{z4}\\
I_4\leq&\sum_{0<|\ell|\leq 3}||[D^{\ell},B\cd] \na c||_{L^2}||\na v||_{H^2}+\sum_{0<|\ell|\leq 3}||[D^{\ell},B\cd] \na v||_{L^2}||\na c||_{H^2}\nonumber\\
\leq&C\Big(||\na B||_{L^\infty}+||\na^3 B||_{L^\infty}\Big)\Big(||v||^2_{H^3}+||c||^2_{H^3}\Big)\label{z5}.
\end{align}
By Leibniz's formula and H\"{o}lder's inequality, one has
\bal
I_5\leq&~||c\cd \na B||_{H^3}||v||_{H^3}+||v\cd \na B||_{H^3}||c||_{H^3}\nonumber\\
\leq&~ C\Big(||\na B||_{L^\infty}+||\na^4 B||_{L^\infty}\Big)\Big(||v||^2_{H^3}+||c||^2_{H^3}\Big),\label{z7}\\
I_6\leq&~||c\cd \na U||_{H^3}||v||_{H^3}+||v\cd \na U||_{H^3}||c||_{H^3}\nonumber\\
\leq&~ C\Big(||\na U||_{L^\infty}+||\na^4 U||_{L^\infty}\Big)\Big(||v||^2_{H^3}+||c||^2_{H^3}\Big)\label{z8}.
\end{align}
Owing to the H\"{o}lder equality, we deduce
\bal
I_{7}\leq&~ \big(||f||_{H^3}+||G||_{H^3}\big)||v||_{H^3}\nonumber\\
\leq&~ C||f,G||_{H^3}+C||f,G||_{H^3}||v||^2_{H^3},\label{z10}\\
I_{8}\leq&~ ||h||_{H^3}||c||_{H^3}\nonumber\\
\leq&~ C||h||_{H^3}+C||h||_{H^3}||v||^2_{H^3}.\label{z11}
\end{align}
Putting all the estimates \eqref{z1}--\eqref{z11} together with \eqref{z0}, we obtain
\bal\label{z12}
&\quad \frac{\dd}{\dd t}\Big(||v||^2_{H^3}+||c||^2_{H^3}\Big)+||\Lambda^\frac{\alpha}{2} v||^2_{H^3}+||\Lambda^\frac{\beta}{2} c||^2_{H^3}\nonumber\\
\les&\Big(||v||_{H^3}+||c||_{H^3}\Big)\Big(||\na v||^2_{H^2}+||\na c||^2_{H^2}\Big)+||f,h,G||_{H^3}\nonumber\\
&\quad +\Big(||B,U,\nabla^4B,\nabla^4U||_{L^\infty}+||f,h,G||_{H^3}\Big)\Big(||v||^2_{H^3}+||c||^2_{H^3}\Big).
\end{align}
Notice that $\mathrm{supp}\ \hat{U_0}\subset \mathcal{C}$ and the fact $(U,B)=(e^{-\mu t}U_0,e^{-\nu t}U_0)$, we can verify
\bal
||\nabla U||_{L^\infty}+||\nabla^4 U||_{L^\infty}\leq& Ce^{-\mu t}||U||_{L^\infty}\leq Ce^{-\mu t}||\hat{U_0}||_{L^1},\label{z13}\\
||\nabla B||_{L^\infty}+||\nabla^4 B||_{L^\infty}\leq& Ce^{-\nu t}||B||_{L^\infty}\leq Ce^{-\nu t}||\hat{U_0}||_{L^1}\label{yyh}.
\end{align}
Inserting \eqref{estimate-l} and \eqref{estimate-2} into \eqref{z12} yields
\bal\label{z15}
&\quad \frac{\dd}{\dd t}\Big(||v||^2_{H^3}+||c||^2_{H^3}\Big)+||\Lambda^\alpha  v||^2_{H^3}+||\Lambda^\beta c||^2_{H^3}\nonumber\\
&\les\Big(||v||^2_{H^3}+||c||^2_{H^3}\Big)^{\fr12}\Big(||\Lambda^\alpha  v||^2_{H^3}+||\Lambda^\beta c||^2_{H^3}\Big)
\nonumber\\&\quad+e^{-\min\{\mu,\nu\}t}\Big(||\hat{U}_0||_{L^1}+\ep||U_0||_{L^2}(1+||\hat{U}_0||_{L^1})\Big)\Big(||v||^2_{H^3}+||c||^2_{H^3}\Big)
\nonumber\\&\quad+e^{-\min\{\mu,\nu\}t}\ep||U_0||_{L^2}(1+||\hat{U}_0||_{L^1}).
\end{align}
Now, we define
\bbal
\Gamma\triangleq\sup\{t\in[0,T^*):\sup_{\tau\in[0,t]}\Big(||v(\tau)||^2_{H^3}+||c(\tau)||^2_{H^3}\Big)\leq \eta\},
\end{align*}
where $\eta$ is a small enough positive constant which will be determined later on.

Assume that $\Gamma<T^*$. For all $t\in[0,\Gamma]$, we obtain from \eqref{z15} that
\bbal
\frac{\dd}{\dd t}\Big(||v||^2_{H^3}+||c||^2_{H^3}\Big) &\leq
Ce^{-\min\{\mu,\nu\}t}\Big(||\hat{U}_0||_{L^1}+\ep||U_0||_{L^2}(1+||\hat{U}_0||_{L^1})\Big)\Big(||v||^2_{H^3}+||c||^2_{H^3}\Big)
\nonumber\\&\quad+Ce^{-\min\{\mu,\nu\}t}\ep||U_0||_{L^2}(1+||\hat{U}_0||_{L^1})\Big),
\end{align*}
which follows from the assumption \eqref{condition} that
\bbal
||v||^2_{H^3}+||c||^2_{H^3}&\leq C(||v_0||^2_{H^3}+||c_0||^2_{H^3}+\ep||U_0||_{L^2}(1+||\hat{U}_0||_{L^1}))
\\&\quad \times\exp\Big( C(||\hat{U}_0||_{L^1}+\ep||U_0||_{L^2}(1+||\hat{U}_0||_{L^1}))\Big)\leq C\delta.
\end{align*}
Choosing $\eta=2C\delta$, thus we can get
\bbal
\sup_{\tau\in[0,t]}\Big(||v(\tau)||^2_{H^3}+||c(\tau)||^2_{H^3}\Big)&\leq \fr\eta2 \quad\mbox{for}\quad t\leq \Gamma.
\end{align*}
So if $\Gamma<T^*$, due to the continuity of the solutions, we can obtain that there exists $0<\epsilon\ll1$ such that
\bbal
\sup_{\tau\in[0,t]}\Big(||v(\tau)||^2_{H^3}+||c(\tau)||^2_{H^3}\Big)&\leq \fr\eta2 \quad\mbox{for}\quad t\leq \Gamma+\epsilon<T^*,
\end{align*}
which is contradiction with the definition of $\Gamma$.

Thus, we can conclude $\Gamma=T^*$ and
\bbal
\sup_{\tau\in[0,t]}\Big(||v(\tau)||^2_{H^3}+||c(\tau)||^2_{H^3}\Big)&\leq C<\infty \quad\mbox{for all}\quad t\in(0,T^*),
\end{align*}
which implies that $T^*=+\infty$. This completes the proof of Theorem \ref{the1.1}. $\Box$

\section*{Acknowledgments} J. Li is supported by the National Natural Science Foundation of China (Grant No.11801090). Y. Yu is supported by the Natural Science Foundation of Anhui Province (No.1908085QA05).

\end{document}